\documentstyle[12pt,psfig]{article}

\reversemarginpar

\makeatletter
\def\@begintheorem#1#2{\it \trivlist \item[\hskip \labelsep{\bf #1\ #2.}]}
\newtheorem{teo}{Theorem}[section]
\newtheorem{rem}[teo]{Remark}
\newtheorem{lem}[teo]{Lemma}
\newtheorem{cor}[teo]{Corollary}

\newtheorem{prop}[teo]{Proposition} 
 
\def\finedim#1{{\hfill\hbox{\enspace\fbox{\ref{#1}}}}\vspace{5pt}}
\def\dim#1{\noindent{\it Proof of} {\hspace{2pt}}\ref{#1}.}

\font\scpicc=cmcsc10
\font\sc=cmcsc10 scaled 1200
\font\titsc=cmcsc10 scaled 1400

\def\compo{\,{\scriptstyle\circ}\,}

\newfont{\Bbb}{msbm10 scaled 1200}
\def\mr{{\hbox{\Bbb R}}}
\def\mc{{\hbox{\Bbb C}}}
\def\mz{{\hbox{\Bbb Z}}}

\def\mh{{\hbox{\Bbb H}}}

\newfont{\Got}{eufm10 scaled 1200}

\def\uu{{\cal U}}

\def\cc{{\cal C}}
\def\TT{{\cal T}}
\def\mm{{\cal M}}

\def\Mbar{\overline M}
\def\hatM{\widehat M}
\def\defo{{\cal D}}
\def\zzero{z^{(0)}}
\def\interior{{\rm int}}
\def\argum{{\rm arg}}
\def\isompiutre{{\rm Isom}^+(\mh^3)}

\def\GL{{\rm GL}}
\def\imunit{{\sqrt{-1}\,}}
\def\struct{{\hbox{\Got h}}}
\def\simstruct{{\hbox{\Got s}}}
\def\zznbd{{\cal U}}
\def\infnbd{{\cal F}}
\def\cont{{\rm C}}
\def\Aff{{\rm Aff}}
\title{Negatively Oriented Ideal Triangulations and a Proof of Thurston's
Hyperbolic Dehn Filling Theorem}

\author{Carlo {\titsc Petronio} \and Joan {\titsc Porti}}

\begin{document}

\maketitle

\noindent{\small{\scpicc Abstract}. We give a complete proof of Thurston's
celebrated hyperbolic Dehn filling theorem, following the ideal triangulation
approach of Thurston and Neumann-Zagier. We avoid to assume that a genuine
ideal triangulation always exists, using only a partially flat
one, obtained by subdividing an Epstein-Penner decomposition. This
forces us to deal with negatively oriented tetrahedra. Our analysis of the set
of hyperbolic Dehn filling coefficients is elementary and self-contained. In
particular, it does not assume smoothness of the complete point in the variety
of deformations.}

\vspace{.2cm}


\vspace{.7cm}

\noindent Thurston's hyperbolic Dehn filling theorem is one of the
greatest achievements in the geometric theory of 3-dimensional manifolds,
and the basis of innumerable results proved over the last twenty years.
Despite these facts, we do not think that a completely satisfactory
written account of the proof exists in the literature, and the aim of this
note is to help filling a gap which could become embarrassing on the long
run. We follow the approach through ideal triangulations, sketched by
Thurston in his notes~\cite{thurston:notes} and later used by Neumann and
Zagier in their beautiful paper~\cite{ne:za}, to prove volume estimates on
the filled manifolds. However, we modify the argument in~\cite{ne:za}
under two relevant respects, which we will explain in detail in this
introduction, after giving the statement of the result itself. We include
both the ordinary and the cone manifold case. 

\begin{teo}\label{main:teo}

Let $M$ be an orientable, non-compact, complete, finite-volume hyperbolic
$3$-manifold. Denote by $\Mbar$ the compact manifold of which $M$ is the
interior, and by $T_1,\dots,T_k$ the tori which constitute $\partial\Mbar$. For
all $i$, choose a basis $\lambda_i,\mu_i$ of $H_1(T_i)$. Denote by $C$ the set
of coprime pairs of integers, together with a symbol $\infty$. For
$c_1,\dots,c_k\in C$ denote by $M_{c_1\cdots c_k}$ the manifold obtained from
$\Mbar$ as follows: if $c_i=\infty$, remove $T_i$; if $c_i=(p_i,q_i)$, glue to
$\Mbar$ along $T_i$ the solid torus $D^2\times S^1$, with the meridian
$S^1\times\{*\}$ being glued to a curve homologous to $p_i\lambda_i+q_i\mu_i$.
Then:
\begin{enumerate}
\item There exists a neighbourhood $\infnbd$ of $(\infty,\dots,\infty)$ in
$C^k$, where $C$ is topologized as a
subset of $S^2=\mr^2\sqcup\{\infty\}$, such that for
$(c_1,\dots,c_k)\in\infnbd$ the manifold $M_{c_1\dots c_k}$ admits
a complete finite-volume hyperbolic structure.
\item Given any $c_1,\dots,c_k\in C$, for small enough positive real numbers
$\vartheta_1,\dots,\vartheta_k$, the manifold $M_{c_1\dots c_k}$ admits the
structure of a complete finite-volume hyperbolic cone manifold,  with cone
locus given by the cores $\{0\}\times S^1$ of the solid tori glued to the
$T_i$'s such that $c_i\neq\infty$, where the cone angle is $\vartheta_i$.
\end{enumerate}
\end{teo}

The first difference of our proof with respect to~\cite{ne:za} is that we
start from a {\em partially flat} ideal triangulation of $M$, namely one
in which some of the tetrahedra degenerate into flat quadrilaterals with
distinct vertices. The existence of such a triangulation easily follows
from a result of Epstein and Penner~\cite{ep:penner}. The argument
in~\cite{ne:za} was based on the assertion that $M$ is itself obtained by
Dehn filling from a hyperbolic manifold which admits a {\em genuine} ideal
triangulation. The reader was addressed to a pre-print of Thurston, later
published as~\cite{thurston:annals}, for the proof of the assertion, but
the result appears to be missing in the printed form of Thurston's paper. 

Some historical explanation about ideal triangulations is in order here. It was
believed for quite some time by several people that the existence of {\em
genuine} ideal triangulations could be proved as an easy consequence of the
result of Epstein and Penner~\cite{ep:penner}. Eventually, this was recognized
to be false, and general existence presently appears to be an open problem (see
for instance~\cite{giapponesini} for sufficient conditions based on the
Epstein-Penner decomposition, and~\cite{snappea} for experimental evidence).
The first named author is responsible, among others, for the spreading of the
erroneous belief that~\cite{ep:penner} implies existence of triangulations. In
particular, the proof presented in~\cite{libro} of Thurston's hyperbolic Dehn
filling theorem is incomplete, because it assumes from the beginning that a
genuine ideal triangulation exists. 

Starting from an ideal triangulation which is partially flat, it becomes
inevitable, when deforming the structure, to deal with negatively oriented
tetrahedra, {\em i.e.}~to consider positive-measure overlapping of the
geometric tetrahedra. The original part of this paper consists of a
careful analysis of this overlapping phenomenon. In particular, we
explicitly show how to associate to a deformed triangulation a hyperbolic
structure on the manifold, and we describe a developing map for this
structure. Since our main motivation was to give a proof of
Theorem~\ref{main:teo}, we have confined our study to ideal triangulations
of the sort which naturally arises when subdividing an Epstein-Penner
decomposition. It is probably possible to extend this study to general
partially flat triangulations, but we believe that the technical details
could be considerably harder (see Section~\ref{deformation:section}). 

The second difference with~\cite{ne:za} in our approach is that we do not
attempt to prove smoothness of the complete point in the deformation space
of the hyperbolic structure. In~\cite{ne:za} the proof of smoothness again
relies on assertions attributed to Thurston, of which no proof (or even
exact statement) is explicitly provided.
Smoothness can actually be proved in the context of
the {\it representation} rather than {\it triangulation} approach to
deformations, see~\cite{kapo}. As mentioned in~\cite{kapo} and sketched
in~\cite{thurston:notes} and~\cite{indiana}, the Dehn filling theorem can
probably be established using the representation approach only, starting
from smoothness near the complete point. However this approach relies on
technical cohomology computations, so we have preferred to stick to
the more elementary and geometric approach through triangulations. 
Thurston actually claims that smoothness can be established also in the
context of triangulations, looking carefully at the equations which define
the space of deformations (personal communication to the first named
author, Berkeley, June 1998).  Being unable to provide the details for
this argument, we have decided not to establish smoothness, but to modify
the proof in~\cite{ne:za} to a possibly singular context.  Our proof that
the set of ``good'' filling parameters indeed covers a neighbourhood of
$(\infty,\dots,\infty)$ becomes somewhat more involved without assuming
smoothness. It uses classical tools from the theory of stratifications and
analytic spaces, which appear to be more suited to a local argument than
tools coming from algebraic geometry, used for instance in~\cite{indiana}.

\vspace{.5cm}

\noindent{\sc Acknowledgements}. We would like to thank the Universities of
Pisa and Toulouse for travel and financial support during the preparation of
this paper. We gratefully thank Riccardo Benedetti and Michel Boileau for many
helpful conversations. In particular, it is a pleasure to acknowledge that the
proof of Proposition~\ref{infinity:covered:prop} emerged from discussions with
Benedetti. The first named author also thanks the Department of Mathematics of
the University of Parma for its friendly hospitality.

\section{Deformation of partially flat
triangulations}\label{deformation:section} 

We describe in this section how to subdivide an Epstein-Penner
decomposition into a partially flat ideal triangulation, and how to associate
to a modified choice of the moduli of the tetrahedra a deformed hyperbolic
structure.

\paragraph{Convex ideal cellularization.} Let us fix for the rest of the paper
a manifold $M$ as in the statement of Theorem~\ref{main:teo}. See for
instance~\cite{libro} or~\cite{rat} for the appropriate definitions, and for
the proof that indeed $M=\interior(\Mbar)$ with
$\partial\Mbar=T_1\sqcup\dots\sqcup T_k$.  It was proved in~\cite{ep:penner}
that there exist convex ideal polyhedra $P_\alpha$, $\alpha=1,\dots,\nu$, in
$\mh^3$ such that $M$ is obtained from their disjoint union via face-pairings.
Each face-pairing will be an isometry $\varphi_i:F_i\to F'_i$ between a
codimension-1 face $F_i$ of some $P_\alpha$ and one such face $F'_i$ of some
other $P_\alpha$ (possibly the same $P_\alpha$, but $F_i\neq F'_i$). Here $i$
ranges between $1$ and half the total number of faces of the $P_\alpha$'s.
Orientability implies that $\varphi_i$ reverses the induced orientation, where
the $P_\alpha$'s are oriented as subsets of $\mh^3$. One way to express the
fact that $M=\bigsqcup P_\alpha/\{\varphi_i\}$ is to say that the quotient of
$\bigsqcup P_\alpha$ under the equivalence relation generated by the
$\varphi_i$ is homeomorphic to $M$, and, modulo this homeomorphism, the
projection into $M$ of the interior of each $P_\alpha$ is an
orientation-preserving isometry. The reason for spelling out this definition is
that later on we will need to deal with less obvious identification spaces.
See~\cite{ep:petronio} for the most general conditions under which a set of
face-pairings on a set of polyhedra defines a manifold or an orbifold.

\begin{figure}
\centerline{\psfig{file=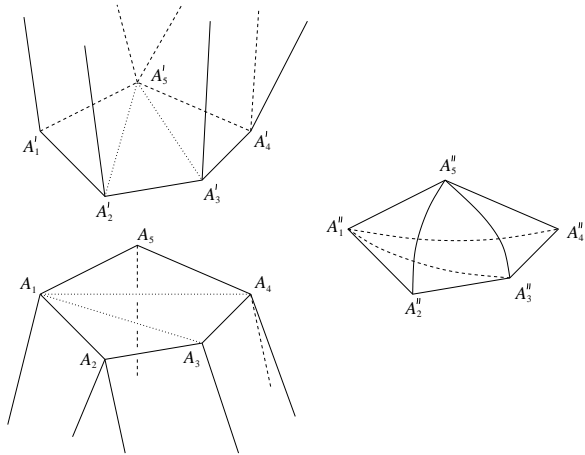,height=6cm}}
\caption{\label{addflat} If two paired pentagonal faces as in the figure 
are subdivided by the dotted lines shown, we add the ``flat'' 
tetrahedra $(A''_1,A''_2,A''_3,A''_5)$ and $(A''_1,A''_3,A''_4,A''_5)$}
\end{figure}

\paragraph{Partially flat triangulation.} We choose now a vertex $v_\alpha$ in
each $P_\alpha$. Moreover, for each of the faces of $P_\alpha$ not containing
$v_\alpha$, we choose a vertex, and take cones from this vertex over the edges
not containing it, to subdivide the face into  triangles. Now we take cones
from $v_\alpha$ over the triangles thus obtained. The result is that $P_\alpha$
has been subdivided into ideal hyperbolic tetrahedra. It will be convenient to
call {\em facets} the triangles into which the original {\it faces} of
$P_\alpha$  have been subdivided. If we now consider a face-pairing
$\varphi_i:F_i\to F'_i$, it may or not be the case that $\varphi_i$ respects
the subdivisions of $F_i$ and $F'_i$ into facets. If subdivisions are not
respected, we can insert geometrically flat ideal tetrahedra between $F_i$ and
$F'_i$, to reconcile these subdivisions, as sketched in  Fig.~\ref{addflat}.

To be precise, assume $F_i$ and $F'_i$ have been triangulated by taking
cones over vertices $w_i$ and $w'_i$ respectively. We identify $\partial
F_i$ to $\partial F'_i$ via $\varphi_i$, and refer to some abstract
version $\gamma_i$ of this loop, disjoint from the original polyhedra. If
$w_i=w'_i$ in $\gamma_i$ then the triangulations of $F_i$ and $F'_i$
match, and there is nothing to do. If $w_i$ and $w'_i$ are the endpoints
of an edge $e$ of $\gamma_i$, as in Fig.~\ref{addflat}, then for every
edge $e'$ of $\gamma_i$ disjoint from $e$ we add the tetrahedron that is
the join of $e$ and $e'$. In the remaining cases the edge between $w_i$
and $w'_i$ is an interior edge of both the triangulations of the faces
$F_i$ and $F'_i$. Then we divide the faces along this edge and apply twice
the previous construction. 

>From the topological point of view, we are led to consider the ideal
triangulation $\TT$ of $M$ which consists of all the ``fat'' tetrahedra
obtained by subdividing the  $P_\alpha$, together with the ``flat'' tetrahedra
just inserted. Recall that a {\em topological ideal triangulation} of $M$ is
just a collection of orientation-reversing  simplicial pairings between the
faces of a finite number of copies of the standard tetrahedron, with the
property that the identification space  defined by the pairings is homeomorphic
to the space $\hatM$ obtained from $\Mbar$ by collapsing each boundary
component to a point. In particular, the name ``fat'' or ``flat'', used for a
tetrahedron of $\TT$, only refers to the way the tetrahedron arose from the 
original geometric subdivision of $M$. The tetrahedron in its own right, as a
member  of $\TT$, is always ``fat''. 

Even if the topological triangulation $\TT$ does depend on the initial choice
of vertices on  the $P_\alpha$, we will fix one such choice and refer to a
definite $\TT$.

\begin{figure}
\centerline{\psfig{file=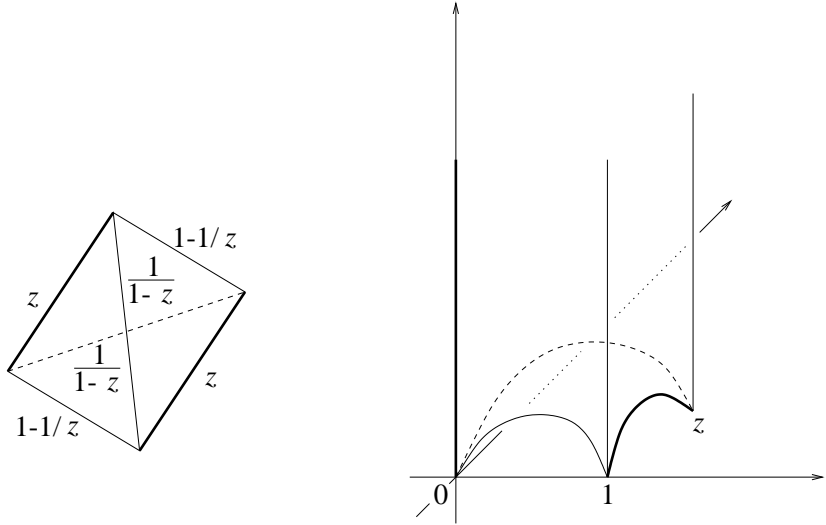,height=7cm}}
\caption{\label{tetramod} Moduli of an ideal tetrahedron, using the 
$\mc\times(0,\infty)$ model of $\mh^3$}
\end{figure}

\paragraph{Consistency and completeness equations.} Recall now that if we fix a
pair of opposite edges on the standard ideal tetrahedron $\Delta$, the
realizations of $\Delta$ as an oriented ideal tetrahedron in $\mh^3$ are
parametrized (up to oriented isometry) by the upper half-plane
$\pi_+=\{z\in\mc:\ \Im(z)>0\}$,  as described in Fig.~\ref{tetramod}. 
This correspondence easily extends to $\mr\setminus\{0,1\}$ to cover the case
where $\Delta$ flattens out to a quadrilateral with distinct vertices. We will
interpret parameters in $-\pi_+$ as describing tetrahedra with negative
orientation (in particular, negative volume).

Given an ideal triangulation $\TT$ of $M$ consisting of tetrahedra
$\Delta_1,\dots,\Delta_n$, we can fix a pair of opposite edges on each
$\Delta_j$, choose a modulus $z_j\in\pi_+$ and ask ourselves if $M$ admits a
(complete) hyperbolic structure inducing on each $\Delta_j$ the structures
described by $z_j$. The answer, which goes back to 
Thurston~\cite{thurston:notes} (see also~\cite{libro}), is given by two systems
of equations in $z=(z_1,\dots,z_n)$. We first have the 
{\em consistency}
equations $\cc^*_\TT(z)$, which prescribe that the  product of the moduli
around each edge should be 1 and the sum of the corresponding arguments should
be $2\pi$. The system $\cc^*_\TT(z)$ is satisfied if and only if there
exists on $M$ a (possibly incomplete) hyperbolic structure as mentioned. In
practice one often needs to consider only the system $\cc_\TT(z)$  obtained by
neglecting the condition on arguments, because close enough to a solution
$\zzero$ of $\cc^*_\TT$, the systems $\cc_\TT$ and $\cc^*_\TT$ are
equivalent. The other equations $\mm_\TT(z)$ one needs to consider, called {\em
completeness} equations, are rational equations in $z$ determined by the
combinatorics of $\TT$, just as it happens for $\cc_\TT$. They have a
geometrical meaning only when $\cc^*_\TT(z)$ holds. In this case a
representation $\rho_z$ of $H_1(\partial\Mbar)$  into the group of affine
automorphisms of $\mc$  is well-defined up to conjugation, and $\mm_\TT(z)$
means that the image of $\rho_z$ consists of translations. 
An exact combinatorial description of $\mm_\TT(z)$ is provided after the
statement of Theorem~\ref{new:main:teo}.

\paragraph{Partially flat and negatively oriented solutions.} The  geometric
meaning of $\cc^*_\TT(z)$ and $\mm_\TT(z)$ for $z_1,\dots,z_n\in\pi_+$ is
as follows. First one realizes the abstract face-pairings as isometries between
the faces of the ideal tetrahedra in $\mh^3$ corresponding to $z_1,\dots,z_n$.
The resulting identification space is homeomorphic to $M$, and a hyperbolic
structure is defined away from the edges. Consistency equations
$\cc^*_\TT(z)$  translate the fact that this structure extends to edges,
and  $\mm_\TT(z)$ translates completeness. As mentioned, the resulting systems
$\cc_\TT$ and $\mm_\TT$ are rational and  depend only on the combinatorics of
$\TT$. Moreover only denominators $z_j$ and $1-z_j$ appear, so it makes sense
to consider solutions $z\in(\mc\setminus\{0,1\})^n$. This is not quite the case
for $\cc^*_\TT(z)$, because for $z\in-\pi_+$ there is no obvious way to
choose arguments for $z,1/(1-z),1-1/z$ so that their sum gives $\pi$. We will
deal with this small subtlety below.

Even if one disregards the problem about arguments, the geometric
interpretation of a solution $z$ of $\cc_\TT$ is not so clear when some $z_j$
is not in $\pi_+$. The idea is that if $z_j\in-\pi_+$ then $\Delta_j$ should
overlap with some other $\Delta_{j'}$ with $z_{j'}\in\pi_+$ (actually, at least
two of them, so that the algebraic number of tetrahedra covering each point is
always 1), but it is not easy to turn this idea into a general formal
definition. Actually, a general definition cannot work, as the following
discussion shows. Consider the tame case where some $z_j$ are in $\pi_+$ and
some (but not all) are in $\mr\setminus\{0,1\}$. If we take the corresponding
``fat'' and ``flat''  tetrahedra in $\mh^3$, we can still glue their faces
together, but it was shown in~\cite{pe:we} that the resulting identification
space is in general not homeomorphic to $M$. If moduli in $-\pi_+$ are
involved, the situation can of course get even worse.

\begin{figure}
\centerline{\psfig{file=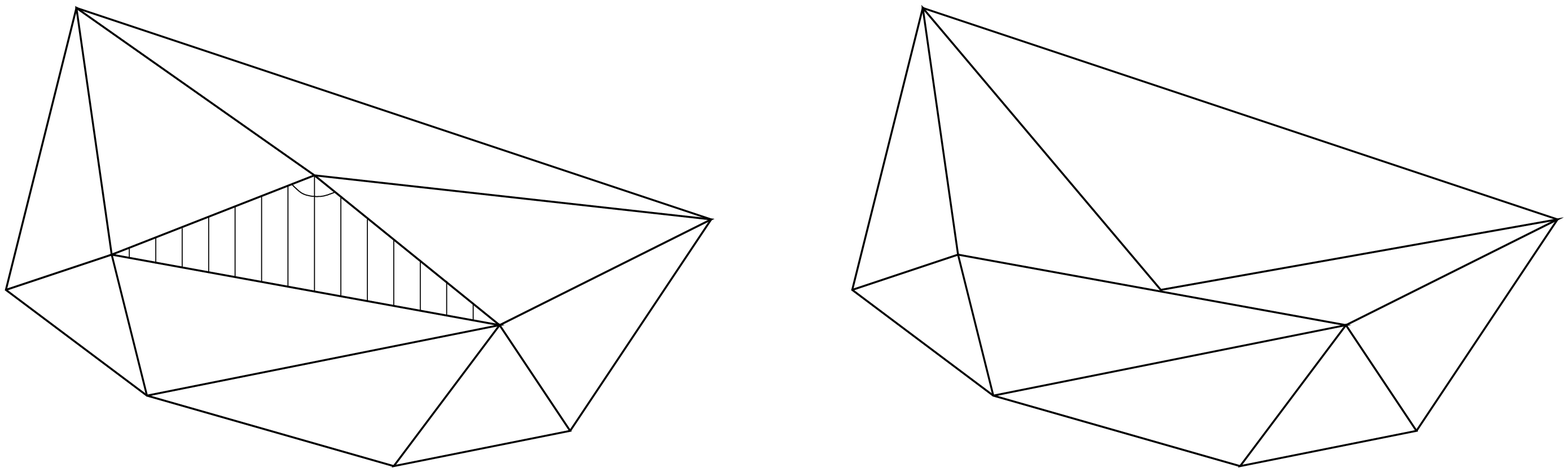,height=2cm}}
\caption{\label{folitria}Foliation representing the flattening of a triangle}
\end{figure}

\paragraph{The complete solution.} We note first that for
$z\in\mr\setminus\{0,1\}$ there is an obvious good choice for the arguments of
$z,1/(1-z),1-1/z$, namely $\argum(t)=\pi$ for $t<0$ and $\argum(t)=0$
otherwise. So, it makes sense to consider partially flat solutions $z$ of
$\cc^*_\TT$. As mentioned, such a $z$ does not have in general a geometric
meaning. However, it was shown in~\cite{pe:we} that if $z$ is a solution also
of $\mm_\TT$ then the identification space obtained from the fat and flat
tetrahedra is  indeed $M$, and a complete hyperbolic structure is naturally
defined. This result itself is not used in this paper, but we will employ the
following technical tool introduced in~\cite{pe:we} for the proof. To signify
the flattening of a genuine triangle into a segment we will foliate the
triangle, as sketched in  Fig.~\ref{folitria}.
One of the main points in~\cite{pe:we} is the proof that the simultaneous
collapse of all the foliated components does not alter the topology.

Going back to the specific situation arising from the subdivision of an
Epstein-Penner decomposition of $M$, we see that we can assign
a modulus $\zzero_j$ to each $\Delta_j$ in $\TT$, where $\zzero_j\in\pi_+$
if $\Delta_j$ lies in some $P_\alpha$, and $\zzero_j\in\mr\setminus\{0,1\}$
if $\Delta_j$ is one of the tetrahedra we have inserted.

\begin{lem}\label{zzero:vera:soluz}
$\zzero$ is a solution of $\cc^*_\TT$ and $\mm_\TT$.
Moreover the foliated components arising on $\partial\Mbar$ have one 
of the shapes described in Fig.~\ref{tamefoli}.
\begin{figure}
\centerline{\psfig{file=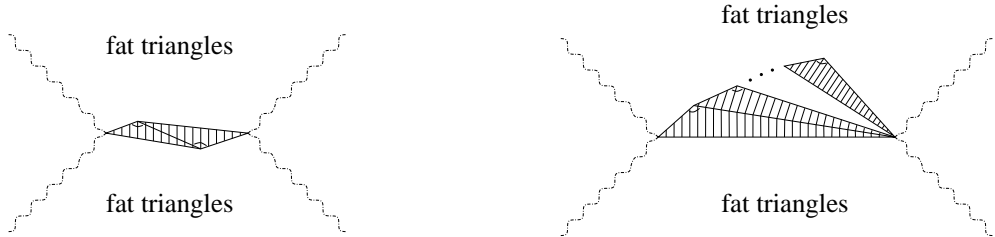,height=3.1cm}}
\caption{\label{tamefoli}Foliated components arising from subdivision of an 
Epstein-Penner decomposition}
\end{figure}
\end{lem}

\dim{zzero:vera:soluz} The first assertion is obvious: we already know that $M$
is complete hyperbolic, and $\zzero$ corresponds to a geometric partially flat
triangulation, so the geometric interpretation of $\cc^*_\TT$ and $\mm_\TT$
is the same as for genuine triangulations. The second assertion is easily
proved by taking transversal sections in Fig.~\ref{addflat} near the ideal
vertices. \finedim{zzero:vera:soluz}

Foliated components as in Fig.~\ref{tamefoli} are called {\em bigons}. 
Lemma~\ref{zzero:vera:soluz} implies that the foliated components on
$\partial\Mbar$ corresponding to $\zzero$ are bigons intersecting each other
only at their ends. This fact will be used in the sequel.

\begin{rem} {\em If one considers a general partially flat solution of
$\cc^*_\TT$ and $\mm_\TT$, annular foliated components  and more
complicated intersections between components can appear on $\partial\Mbar$,
see~\cite{pe:we}. This makes the analysis of the deformed structures
considerably harder, and explains why we have decided to concentrate on
solutions arising from Epstein-Penner decompositions.} \end{rem}

\paragraph{Solutions near the complete solution.} From now on we will only be
concerned with solutions $z$ of $\cc_\TT$ lying in an arbitrarily small
neighbourhood $\zznbd$ of $\zzero$. Formally, all our statements should contain
the phrase ``$\zznbd$ can be taken so small that...'', but we will omit it
systematically. We define $\defo=\{z\in\zznbd:\ \cc_\TT(z)\}$. We note first
that on $\zznbd$ the arguments  can be defined by continuity also for the
moduli in $-\pi_+$, and of course the  resulting system $\cc^*_\TT$ is
equivalent to $\cc_\TT$. For this reason we will henceforth leave the
discussion of arguments in the background. Moreover we will assume that for
$z\in\zznbd$, if $\zzero_j\in\pi_+$, then also $z\in\pi_+$. In other words,
flat tetrahedra can become fat, flat, or negative, but fat tetrahedra stay fat.

It will be convenient to denote the generic abstract element of $\TT$ by
$\Delta_j$, and by $\Delta_j(z)$ the geometric version of $\Delta_j$
corresponding to $z\in\defo$. As mentioned, for $z_j\in-\pi_+$ one imagines
$\Delta_j(z)$ to be negatively oriented, but  we will only use directly those
$\Delta_j(z)$ for which $\zzero_j$, and hence  $z_j$, lies in $\pi_+$. We will
also use $P_\alpha(\zzero)$ to emphasize that we are considering the geometric
polyhedron rather than the abstract one $P_\alpha$. For all $\alpha$, let 
$J_\alpha$ be the set of indices $j$ such that $\Delta_j$
appears in the original subdivision of $P_\alpha$. 
Consider also the set of face-pairings $p_\alpha$ corresponding to the
triangles lying in the interior of $P_\alpha$. In this context a face-pairing
is just a combinatorial rule,  but when the abstract tetrahedra are turned into
geometric ones, an isometry is uniquely  determined.

\begin{lem}\label{non:convex:polyhedra} For $z\in\defo$ and for all $\alpha$,
the tetrahedra $\Delta_j(z)$, $j\in J_\alpha$ can be assembled along $p_\alpha$
to give a (probably non-convex) ideal polyhedron $P_\alpha(z)$ in $\mh^3$ with
triangular faces, combinatorially equivalent (in particular, homeomorphic) to 
$P_\alpha$ (endowed with the {\em facets} structure). \end{lem}

\dim{non:convex:polyhedra} We first note that, using the projective model of
$\mh^3$, ideal polyhedra can be viewed as compact Euclidean polyhedra with
vertices on the unit sphere.  Choosing a maximal tree in the graph
corresponding to the pairing $p_\alpha$, we can realize in $\mh^3$ the
$\Delta_j(z)$, $j\in J_\alpha$, so that the pairings in the tree are given by
actual overlapping. Moreover we can define a map 
$f_j(z):\Delta_j(\zzero)\to\Delta_j(z)$, for instance using Euclidean
coordinates and taking convex combinations of vertices. Since $z$ satisfies
$\cc_\TT(z)$, these $f_j(z)$'s match to give a map $F_\alpha(z):P_\alpha(\zzero)\to
P_\alpha(z)$. Moreover $F_\alpha(z)$ is locally injective. To conclude we note
that $P_\alpha(z)$  converges to the identity of $P_\alpha(\zzero)$ as $z$ goes
to $\zzero$, and we use Euclidean  compactness of $P_\alpha(\zzero)$ to deduce
that $F_\alpha(z)$ is eventually injective. All conclusions easily follow.
\finedim{non:convex:polyhedra}

Using the combinatorial equivalence  between $P_\alpha(z)$ and $P_\alpha$, we
can define the faces $F_i(z)$ also for the $P_\alpha(z)$. Each  $F_i(z)$ will
be a (probably non-planar) union of facets.

We define now an {\em abstract} polyhedron $\tilde P_\alpha$ by adding to
$P_\alpha$ all the flat tetrahedra $\Delta_j$ arising from faces $F_i$
contained in $\partial P_\alpha$. Recall that we have artificially broken the
symmetry of face-pairings using the notation $F_i,F'_i$ for a pair of faces to
be glued, so each flat tetrahedron is used once. The $\tilde P_\alpha$ have a
natural facet structure on their boundary. Moreover, using the pairing of
triangles in $\TT$, we deduce a pairing of the facets of the $\tilde P_\alpha$,
and the result of all these facet-pairings is $M$.

The idea is now to replace each $P_\alpha(z)$ by some $\tilde P_\alpha(z)$
having the same combinatorial structure as $\tilde P_\alpha$, so to obtain
$M$ from geometric polyhedra. As obvious, $\tilde P_\alpha(z)$ will result
from elementary modifications on $P_\alpha(z)$, each modification coming
from one of the faces $F_i$ contained in $\partial P_\alpha$. The
elementary modification is itself obvious: $\tilde P_\alpha(z)$ will have the
same vertices on $\partial\mh^3$ as $P_\alpha(z)$, but facets (convex
envelopes of triples of these vertices) will be taken according to the
combinatorial structure of $\tilde P_\alpha$ rather than $P_\alpha$. For
example, consider the situation of Fig.~\ref{addflat}.  Let $P_\alpha$ be
the polyhedron shown below in the figure, and let $F_i=(A_1,\dots,A_5)$.
The collection of facets of $P_\alpha(z)$ contains the triangles
$(A_1(z),A_2(z),A_3(z))$, $(A_1(z),A_3(z),A_4(z))$,
$(A_1(z),A_4(z),A_5(z))$. Now we replace these triangles by
$(A_5(z),A_1(z),A_2(z))$, $(A_5(z),A_2(z),A_3(z))$,
$(A_5(z),A_3(z),A_4(z))$, leaving all other facets of $P_\alpha(z)$
unchanged. The resulting collection of triangles still bounds an ideal
polyhedron in $\mh^3$, which we take as $\tilde P_\alpha(z)$. We will also
denote by $\tilde F_i(z)$ the union of the modified facets. 

\begin{rem}\label{order:remark}{\em Assume that under a face-pairing 
$\varphi_i:F_i\to F'_i$ no edge of the subdivisions of $F_i$ and $F'_i$ is
matched (as in Fig.~\ref{addflat}). Then the flat tetrahedra inserted come in a
natural order starting from $F_i$ and proceeding towards $F'_i$ (in
Fig.~\ref{addflat}, first $(A''_1,A''_5,A''_4,A''_3)$ and then
$(A''_1,A''_5,A''_3,A''_2)$). The transformation of $P_\alpha(z)$ into $\tilde
P_\alpha(z)$ can be viewed as the result of successive transformations
corresponding to the individual flat tetrahedra. Each transformation consists
in replacing a quadrilateral, bent along one diagonal, with the quadrilateral
having the same perimeter and bent along the other diagonal. If the dihedral
angle at the first diagonal is more than $\pi$ then the modulus of the
corresponding tetrahedron is in $\pi_+$, and the tetrahedron is being added to
$P_\alpha$. If the angle is less than $\pi$, then the modulus is in $-\pi_+$,
and the tetrahedron is being deleted. If the angle is $\pi$, the modulus is in
$\mr\setminus\{0,1\}$ and we are only changing the combinatorial structure of
the facets of $P_\alpha$. When the pairing $\varphi_i$ matches an edge of the
subdivisions, this description must be repeated for both of the polygons into
which $F_i$ is divided by the matching edge.}\end{rem}

\begin{teo}\label{glue:tilde:alpha}
\begin{enumerate}
\item\label{simult:point} The above-described modification of $P_\alpha(z)$ 
can be carried out simultaneously for all faces $F_i$. 
\item\label{comb:equiv:point} The resulting collection $\tilde P_\alpha(z)$,
with the face structure given by the $\tilde F_i(z)$ and the $F'_i(z)$, is
combinatorially equivalent to the original collection $P_\alpha$. 
\item\label{isom:pair:point} Each pairing $\tilde F_i(z)\to F'_i(z)$ can be
realized by an isometry. 
\item\label{structure:point} The identification space resulting from the
pairings is homeomorphic to $M$, and it can be endowed with a hyperbolic
structure compatible with the structure defined on the interior of each $\tilde
P_\alpha(z)$. 
\end{enumerate} 
\end{teo}

\dim{glue:tilde:alpha}  It is again useful to identify hyperbolic ideal
polyhedra with compact Euclidean polyhedra with vertices on the sphere. Using
this point of view, let us consider the 1-skeleton $\Gamma_\alpha(\zzero)$ of a
certain $P_\alpha(\zzero)$.  On $\Gamma_\alpha(\zzero)$ we have certain simple
circuits which correspond to the faces of $P_\alpha$. Note that each circuit is
contained in a plane, and the various planes form dihedral angles strictly less
than $\pi$ at the edges of $\Gamma_\alpha(\zzero)$. Now we consider the same
circuits in the modified 1-skeleton $\Gamma_\alpha(z)$. By compactness, we
easily see that for $z$ close enough to $\zzero$, the convex envelopes of any
two distinct circuits meet at most in a common edge or vertex of
$\Gamma_\alpha(z)$. This shows
points~\ref{simult:point},~\ref{comb:equiv:point} and the first assertion
in~\ref{structure:point}.

We show point~\ref{isom:pair:point} in the special case of Fig.~\ref{addflat},
leaving to the reader the general case. The idea is to somehow realize in
$\mh^3$ the flat tetrahedra. Let $x$ and $y$ be the moduli along the edge
$(A''_1,A''_5)$ of the tetrahedra $(A''_1,A''_2,A''_3,A''_5)$ and
$(A''_1,A''_3,A''_4,A''_5)$ respectively. Note that
$x(\zzero),y(\zzero)\in(1,\infty)$. Now in the half-plane model of $\mh^3$ we
choose $A''_1(z)=\infty$, $A''_5(z)=0$, $A''_4(z)=1$, $A''_3(z)=x(z)$ and
$A''_2(z)=y(z)\cdot x(z)$. Consistency of $z$ along 
$(A_1,A_4)$ and $(A_1,A_3)$ implies that the unique $f\in\isompiutre$
such that $f(A_1(z))=A''_1(z)$, $f(A_5(z))=A''_5(z)$, and
$f(A_4(z))=A''_4(z)$, also enjoys $f(A_3(z))=A''_3(z)$ and
$f(A_2(z))=A''_2(z)$. Similarly
consistency along $(A_5,A_2)$ and $(A_5,A_3)$ implies that  
$g(A'_l(z))=A''_l(z)$, $l=1,\dots,5$,
for some $g\in\isompiutre$.
Now, the description of $\tilde P_\alpha$ given
in Remark~\ref{order:remark} implies that 
$$\tilde F_i(z)=f^{-1}\Big((A''_5(z),A''_1(z),A''_2(z))\cup
(A''_5(z),A''_2(z),A''_3(z))\cup(A''_5(z),A''_3(z),A''_4(z))\Big)$$
whence the conclusion.

The second assertion in point~\ref{structure:point} follows from
point~\ref{isom:pair:point} and consistency along the original edges of the
$P_\alpha$.\finedim{glue:tilde:alpha}

\section{Developing map and completion\\ of deformed
structures}\label{developing:section} 

We will denote in the sequel by $\struct(z)$ the hyperbolic structure on $M$
constructed in Theorem~\ref{glue:tilde:alpha} for $z\in\defo$. In this section
we will analyze the completion of $\struct(z)$, the key ingredient being the
understanding of the developing map of cusps. We will first  give the statement
needed in Section~\ref{parameters:section} to conclude the proof of
Theorem~\ref{main:teo}, then we will switch to a 2-dimensional setting, and
later we will use the 2-dimensional construction to understand $\struct(z)$. 

\paragraph{Statements of results.} Let us
return to the notation of Theorem~\ref{main:teo} 
and slightly modify it so to unify the two assertions. Consider the set
$$G=\{\infty\}\cup\{g\in\mr^2:\ g=r\cdot(p,q){\rm\ for\ some\ }
r>0{\rm\ and\ relatively\ prime\ }p,q\in\mz\}.$$
(The motivation for the notation is that $G$ consists of {\em G}eneralized
filling coefficients, as opposed to the genuine {\em C}oefficients of the
set $C$ defined in Theorem~\ref{main:teo}.) For $g\in G\setminus\{\infty\}$
note that its expression as  $r\cdot(p,q)$ is unique, and define $c(g)=(p,q)$,
$\vartheta(g)=2\pi/r$. Set $c(\infty)=\infty$. Topologize $G$ as a subset of
$\mr^2\cup\{\infty\}=S^2$. We can now restate Theorem~\ref{main:teo} as
follows:

\begin{teo}\label{new:main:teo}
Under the assumptions of Theorem~\ref{main:teo} there exists a neighbourhood
$\infnbd$ of $(\infty,\dots,\infty)$ in $G^k$ such that for
$(g_1,\dots,g_k)\in\infnbd$ the manifold $M_{c(g_1)\dots c(g_k)}$ admits the
structure of a complete finite-volume hyperbolic cone manifold,  with cone
locus given by the cores $\{0\}\times S^1$ of the solid tori glued to the
$T_i$'s such that $g_i\neq\infty$, where the cone angle is $\vartheta(g_i)$.
\end{teo}

It is perhaps worth noticing here that this statement is actually
independent of the choice of the basis $\lambda_i,\mu_i$ of
$H_1(T_i)$. In fact, a different choice is related through a matrix in
$\GL(2,\mz)$, which induces a homeomorphism of $S^2$ and preserves
coprimality of integer pairs, and hence the function $\vartheta:G\to\mr_+$
introduced above.

Theorem~\ref{new:main:teo} 
is the result which we will establish in the rest of the paper. To
summarize the content of the present section, we now go back to the notation of
Section~\ref{deformation:section}. Note first that for $z\in\defo$ a
homomorphism $h_i(z):H_1(T_i)\to\mc^*$ is defined by
$h_i(z)([\gamma])=(-1)^{\#\gamma_0}L_z(\gamma)$, where $\gamma$ is a simplicial
loop with respect to the triangulation of $T_i$ induced by $\TT$, $\#\gamma_0$
is the number of vertices of $\gamma$ and $L_z(\gamma)$ is the product of all
moduli along angles which $\gamma$ leaves on its left on $T_i$. Recall that $\mm_\TT(z)$ is the
system $\{h_i(z)(\lambda_i)=h_i(z)(\mu_i)=1,\ i=1,\dots,k\}$. Note that
$h_i(\zzero)(\lambda_i)=h_i(\zzero)(\mu_i)=1$, so we can use the holomorphic
branch $\log$ of the logarithm function enjoying $\log(1)=0$ to define maps
$u_i,v_i:\defo\to\mc$ as $u_i(z)=\log(h_i(z)(\lambda_i))$ and
$v_i(z)=\log(h_i(z)(\mu_i))$. We will establish the following:

\begin{teo}\label{completion:teo}\begin{enumerate} 
\item\label{symmetric:u:v:point} For $z\in\defo$, we have $u_i(z)=0$ if and
only if $v_i(z)=0$. 
\item\label{rigidity:u:point} If $z\in\defo$ and $u_1(z)=\cdots=u_k(z)=0$ then
$z=\zzero$.
\item\label{non:real:limit:point} The following limit exists and is not real:
$$\tau_i=\lim_{z\in\defo,u_i(z)\neq 0,z\to\zzero}\quad 
{v_i(z)\over u_i(z)}.$$
\item\label{completion:description:point} Fix $z\in\defo$, and let
$g_1,\dots,g_k\in G$ be such that $g_i=\infty$ when $u_i(z)=0$,
and $g_i=(p_i,q_i)$ with $p_i\cdot u_i(z)+q_i\cdot v_i(z)=2\pi\imunit$
otherwise. Then the completion of $M$ with respect to $\struct(z)$ is
homeomorphic to $M_{c(g_1)\dots c(g_k)}$, and the structure of $M$
extends to a hyperbolic cone manifold structure as described in
Theorem~\ref{new:main:teo}. 
\end{enumerate}\end{teo}

\paragraph{Partially flat triangulations of the torus.} Let us consider a
triangulation $\TT$ of the torus $T$ (all notation overlaps between this
paragraph and the previous section are intentional, and their motivation should
be clear to the reader). The combinatorics of $\TT$ allows to write down
systems $\cc^*_\TT$ and $\mm_\TT$, the latter requiring the choice of a
basis $\lambda,\mu$ of $H_1(T)$. For $z_1,\dots,z_n\in\pi_+$, $\cc_\TT(z)$
holds if and only if there is on $T$ a similarity structure inducing on the
$j$-th triangle the structure with modulus $z_j$. Moreover, also $\mm_\TT(z)$
holds if and only if this structure is compatible with a Euclidean structure.
Let us fix now a solution $\zzero$ of $\cc^*_\TT$ and $\mm_\TT$ which is
only {\em partially} (but not totally) {\em flat}. It was shown in~\cite{pe:we}
that $\zzero$ still yields a Euclidean structure (up to scaling) on $T$.
However, being only interested in the situations arising on $\partial\Mbar$
when subdividing an Epstein-Penner decomposition, we may take as an {\em
assumption} that there is on $T$ a Euclidean structure inducing on the $j$-th
triangle of $\TT$ the structure with modulus $\zzero_j$. Of course when
$\zzero_j$ is real this means that the triangle has been collapsed to a
segment. We will use foliations to signify collapse, as in Fig.~\ref{folitria}.
We will also assume that foliated components of $T$ are bigons intersecting at
their ends only, as in Fig.~\ref{tamefoli}.

Before proceeding, we need to recall that for a solution
$z\in(\mc\setminus\{0,1\})^n$ of $\cc^*_\TT$, a representation
$h(z):H_1(T)\to\mc^*$ can be defined as explained above. Moreover
$\mm_\TT(z)$ is the system $h(z)(\lambda)=h(z)(\mu)=1$. (It follows from
this that all systems $\mm_\TT(z)$ arising from different choices of the
basis of $H_1(T)$ are equivalent to each other. However, we will not
need to change basis.)

\begin{prop}\label{develop:torus}
There 
exist a neighbourhood $\zznbd$ of $\zzero$ in $(\mc\setminus\{0,1\})^n$
such that:
\begin{enumerate}
\item\label{torus:forget:tilde:point} For $z\in\zznbd$, $\cc_\TT(z)$ is
equivalent to $\cc^*_\TT(z)$.
\item\label{torus:lambda:mu:same:point} 
If $\defo:=\{z\in\zznbd: \cc_\TT(z)\}$
and $z\in\defo$, then $h(z)(\lambda)=1$ if and only if $h(z)(\mu)=1$.
\item\label{torus:limit:point} If $u(z)=\log(h(z)(\lambda))$ and
$v(z)=\log(h(z)(\mu))$, where $\log$ is holomorphic near $1\in\mc$ and
$\log(1)=0$, then the limit of $v(z)/u(z)$, as $z$ tends to $\zzero$ in 
$\defo$
and $u(z)\neq 0$, exists and is a non-real number $\tau$.
\item\label{torus:similarity:point} Each $z\in\defo$ defines on $T$ a
similarity structure $\simstruct(z)$.
\item\label{torus:Euclidean:point} For $z\in\defo$, $\simstruct(z)$ is
compatible with a Euclidean structure on $T$ if and only if
$h(z)(\lambda)=1$.
\item\label{torus:developing:point} If $h(z)(\lambda)\neq 1$,
 a developing map
for $\simstruct(z)$ is given by
$$\mr^2\ni(x,y)\mapsto \exp(u(z) x+ v(z) y)\in\mc$$
where $\mr^2$ is the universal cover of $T$, with deck 
transformation group
$\mz^2$. 
\end{enumerate}
\end{prop}

\begin{rem}{\em Oriented similarity structures in dimension two 
 are equivalent to complex affine structures
in dimension one, and we will use both indistinctly.}\end{rem}

\dim{develop:torus} Point~\ref{torus:forget:tilde:point} is clear from
continuity.  To prove the other points we start with the Euclidean structure on
$T$. Let $\sigma_1,\ldots,\sigma_r$ be the triangles of $\TT$ that have
non-zero area for this Euclidean structure. The remaining triangles are
 flat,
so 
they have a longest edge (the one with angle zero at each endpoint), and we
abstractly glue each one of these triangles to its neighbour along the longest
edge.
Since we assume that the foliated components of $T$ are bigons as
described in Fig.~4, each flat triangle is glued to either a fat one or to a
family of flat triangles glued to a fat one. The result of this gluing process
is a family of abstract triangulated polygons $\tilde\sigma_1,\ldots,\tilde\sigma_r$,
such that each 
$\tilde\sigma_i$ contains exactly one
triangle that is fat for the Euclidean structure.
Moreover we have a family of pairings
between the edges of the $\tilde\sigma_i$'s, yielding $T$ as identification
space.

The parameters $z\in\zznbd$ define a complex affine structure on the triangles
$\sigma_i$ that we denote by $\sigma_i(z)$. Now we define the induced
structures on the $\tilde\sigma_i$'s. We first repeat 
  geometrically the combinatorial construction of $\tilde\sigma_i$,
namely we add the triangles with parameter in $\pi_+$ and we remove the
triangles with parameter in $-\pi_+$. The triangles with real parameter are the
flat ones, and for them we add a new
vertex in the interior of the edge they represent, according to the real
parameter. This process is only possible when $z$ is close to $\zzero$.  
We denote
by $\tilde\sigma_i(z)$ the complex affine polygon obtained in this way. The
next lemma proves point~\ref{torus:similarity:point} of
Proposition~\ref{develop:torus}.

\begin{lem}\label{complex:affine:triangles} For $z\in\defo$ and for
$i=1,\ldots,r$, $\tilde\sigma_i(z)$ defines a complex affine structure on the
polygon $\tilde\sigma_i$. These structures match  under the
edge-pairings and induce a complex
affine structure $\simstruct(z)$ on $T$. \end{lem}

\dim{complex:affine:triangles} For the first assertion we have to show that
there is a natural combinatorial equivalence between $\tilde\sigma_i$ and
$\tilde\sigma_i(z)$.  We view
$\partial\tilde\sigma_i(\zzero)$ not as a triangle but as a polygon
combinatorially equivalent to $\partial\tilde\sigma_i$, because each time we
glue a flat triangle we are adding a new vertex. Hence 
$\tilde\sigma_i(\zzero)$ is a polygon in $\mc$,  with every angle but three
equal to $\pi$. Now the polygon $\partial\tilde\sigma_i(\zzero)$  is
combinatorially isomorphic to the abstract polygon
$\partial\tilde\sigma_i$, and $\partial\tilde\sigma_i(z)$ is  isomorphic to
$\partial\tilde\sigma_i(\zzero)$  for $z\in\zznbd$, because the vertices depend
continuously on $z$, so $\partial\tilde\sigma_i(z)$ is equivalent to
$\partial\tilde\sigma_i$. 

Having shown that the
$\tilde\sigma_i(z)$'s are equivalent to the $\tilde\sigma_i$'s, we can now
realize the edge-pairings by similarities. Consistency equations
$\cc_\TT(z)$
are readily seen to imply that the similarity structure defined on $T$ minus
the vertices extends to the vertices, whence the
conclusion.\finedim{complex:affine:triangles}

We next consider the holonomy of $\simstruct(z)$. 
This is a homomorphism $\pi_1(T)\to\Aff(\mc)$
well-defined up to conjugation. Speaking of holonomy we need to refer to
$\pi_1(T)$, but we will freely use the canonical isomorphism with $H_1(T)$.
Given $f\in\Aff(\mc)$, if $f(w)=\alpha w+\beta$ we call $\alpha$ the linear
part of $f$. Note that $\alpha$ is invariant under conjugation, so the
linear part of the holonomy is a well-defined homomorphism $\pi_1(T)\to\mc^*$,
which depends only on the complex affine structure.

\begin{lem}\label{holonomy}
Given $g\in\pi_1(T)$, the linear part of the holonomy of $g$ corresponding to
$\simstruct(z)$ is $h(z)(g)$, where $h$ is defined as above.
In addition, there exists a representative $\rho(z)$ of the holonomy such that
$\rho(z)(\lambda)$ and $\rho(z)(\mu)$ are respectively given by
$$w\mapsto {\rm e}^{u(z)}w+a(z) 
\qquad\hbox{ and } \qquad 
w\mapsto {\rm e}^{v(z)}w+b(z),$$
where $a,b:\defo\to\mc$ are restrictions to $\defo$ of global rational functions
with denominators not vanishing in $\defo$.
\end{lem}

\dim{holonomy} To prove the first assertion we recall the general recipe to
compute the linear part of the holonomy. We consider the CW-decomposition of
$T$ where the 2-cells  are the polygons $\tilde\sigma_1,\ldots,\tilde\sigma_r$
and the 1-skeleton is the union of the boundaries of these polygons. Given an
element in $\pi_1(T)$, we represent it by a path $\gamma$ in the 1-skeleton.
Since the 2-cells are polygons of $\mc$ defined up to similarity, the complex
ratio between two 1-cells with a common end is an invariant of the similarity
structure. The linear part of the holonomy of the oriented path $\gamma$ is the
product of the ratios between each pair of consecutive 1-cells of $\gamma$,
taking care of the orientations.  In our
situation, each ratio between consecutive 1-cells is a product of parameters
$z_j$, $1-1/z_j$ or $1/(1-z_j)$, and one can easily check that the linear part
of the holonomy of $\gamma$ is precisely $h(\gamma)(z)$ as defined above.

To prove the second assertion, we fix a polygon $\tilde\sigma_i$ and one of its
edges. We normalize the developing map $D(z)$ so that it maps this edge to the
segment $[0,1]$ in $\mc$. In the fundamental group, we choose the basepoint to be the initial point of
the edge we have fixed, and consider the holonomy  $\rho(z)$ corresponding to
$D(z)$. The first assertion of the lemma and the definition of $u$ imply that
the linear part of $\rho(z)(\lambda)$ is indeed ${\rm e}^{u(z)}$. Moreover, by
our choices, $a(z)=D(z)(\tilde\lambda(1))$, where $\tilde\lambda$ is a lift of
$\lambda$ to the universal covering such that
$\tilde\lambda(0)\in D(z)^{-1}(0)$. Now, $D(z)$ is constructed by patching
together in $\mc$ triangles with moduli $z_j$, with one triangle having
vertices $0$ and $1$. All resulting vertices, in particular
$a(z)=D(z)(\tilde\lambda(1))$, are therefore polynomials in the $z_j$,
$1-1/z_j$ and $1/(1-z_j)$. This implies the conclusion for $\rho(z)(\lambda)$,
and the same argument applies to $\mu$. \finedim{holonomy}

Since the complex affine structure $\simstruct(\zzero)$ is compatible with
a Euclidean structure, $u(\zzero)=v(\zzero)=0$ and $\langle
a(\zzero),b(\zzero)\rangle$ is a lattice in $\mc$. In particular
$a(\zzero),b(\zzero)\in\mc\setminus\{0\}$ and
$\tau=b(\zzero)/a(\zzero)\in\mc\setminus\mr$. Moreover it follows from the
commutativity between $\lambda$ and $\mu$ that: 
$$a(z) ({\rm e}^{v(z)}-1)=b(z) ({\rm e}^{u(z)}-1).$$
Points~\ref{torus:lambda:mu:same:point},~\ref{torus:limit:point}
and~\ref{torus:Euclidean:point} in Proposition~\ref{develop:torus} follow
directly from this equality. We are left to prove
point~\ref{torus:developing:point}. We start with $\simstruct(\zzero)$. 
Since this structure is compatible with a Euclidean one, it is complete,
because every Riemannian structure on a compact manifold is complete. This
means that the structure $\simstruct(\zzero)$ is realized by the quotient
$\mc/\Gamma$, where the lattice $\Gamma=\langle a(\zzero),b(\zzero)\rangle
<\mc$ is the image of $\pi_1(T)$ under the holonomy. Using the fact that
the isotopy class of a homeomorphism of the torus is determined by its
action on the fundamental group, it follows that a developing map for
$\simstruct(\zzero)$ is given by any equivariant homeomorphism between
$\mr^2$ and $\mc$.  Hence a developing map for $\simstruct(\zzero)$,
normalized as in the proof of the previous lemma, is given by
$$\mr^2\ni(x,y)\mapsto
a(\zzero) x+ b (\zzero) y\in\mc.$$ 
By~\cite{weil} or \S 1.7 of~\cite{CEG}, for $z\in\defo$, to give a
developing map of $\simstruct(z)$ it suffices to deform the developing map
of $\simstruct(\zzero)$ to a local embedding equivariant with the
holonomy. The following family of maps has the required properties:
\begin{equation}
\mr^2\ni(x,y)\mapsto\cases{ a(z)\cdot{ \exp((u(z)x+v(z)y))-1\over
\exp(u(z))-1}& if $u(z)\neq 0$\cr a(z) x+ b(z) y &
otherwise}\quad\in\mc.\label{developing:torus} \end{equation} More
precisely, this is a family of maps from $\mr^2$ to $\mc$ parametrized by
$z\in\defo$. This family depends continuously on the parameter
$z\in\defo$, in the sense that if we have a convergent sequence in
$\defo$, then the corresponding sequence of maps converges uniformly on
compact subsets of $\mr^2$ for the $\cont^1$ topology. In addition, the
map corresponding to $z\in\defo$ is equivariant with the holonomy of
$\simstruct(z)$ in Lemma~\ref{holonomy}. Hence it is a developing map for
$\simstruct(z)$ when $z\in\defo$. When $u(z)\neq0$, if we compose the map
in~(\ref{developing:torus})  with a suitable complex affine
transformation, we obtain the map in
point~\ref{torus:developing:point} of the proposition.
\finedim{develop:torus}

\paragraph{3-dimensional developing map.} Points~\ref{symmetric:u:v:point}
and~\ref{non:real:limit:point} of Theorem~\ref{completion:teo} follow directly
from Proposition~\ref{develop:torus}, considering the various $T_i$'s. To
establish the other points, we go back now to the setting of
Section~\ref{deformation:section}. We know that each $z\in\defo$ defines on $M$
a hyperbolic structure $\struct(z)$, and our plan here is to develop it to
analyze its completion. We will cut $M$ along a collection of disjoint
boundary-parallel tori, getting a compact manifold $M_0$ with boundary,
together with cusps $C_1,\dots,C_k$, with $C_i\cong T_i\times[0,\infty)$ and
$T_i$ corresponding to $T_i\times\{\infty\}$. We will allow ourself to isotope
the cutting tori without changing notation.   Since $M$ is
$\partial$-incompressible, if we take a developing map of $M$ and restrict it
to a component of the preimage (under the universal covering) of $C_i$, we get
a developing map for the restriction $\struct_i(z)$ of $\struct(z)$ to $C_i$.
Therefore the completion of $M$ is obtained by completing the various $C_i$'s
separately and then glueing back to $M_0$ along the tori.

\begin{prop}\label{develop:cusp} If $z\in\defo$ then $\struct_i(z)$ is 
complete
if and only if $u_i(z)=0$. 
 If
$u_i(z)\neq 0$ then
a developing map for $\struct_i(z)$ is given by:
$$
\begin{array}{rcl} \mr^2\times[0,\infty)&\to&\mh^3\cong\mc\times(0,\infty)\\
(x,y,t)&\mapsto &(\exp(u(z) x+ v(z) y),\exp(t+\Re( u(z) x+ v(z) y))).
\end{array} $$ If $p_i\cdot u_i(z)+q_i\cdot v_i(z)=2\pi\imunit/r_i$ for some
coprime pair of integers $(p_i,q_i)$ and a real number $r_i>0$, then the
completion of $C_i$ is obtained by attaching $D^2\times S^1$ to
$T_i\times[0,\infty]$ along $T_i\times\{\infty\}$, with $S^1\times\{*\}$ 
glued
to $(p_i\lambda_i+q_i\mu_i)\times\{\infty\}$, and the result has the
 structure
of a hyperbolic cone manifold with boundary, with cone locus
 $\{0\}\times S^1$ and angle $2\pi/r_i$.\end{prop}

\dim{develop:cusp} We will use both the statement and the proof of
Proposition~\ref{develop:torus}, denoting by $\simstruct_i(z)$ the 
similarity structure defined on $T_i$ according to that proposition.
Now, the hyperbolic
structure $\struct(z)$ on the open manifold $M$ induces another
similarity structure on
$T_i$, which we denote by $\simstruct^*_i(z)$. We have the
following:

\begin{lem}\label{same:structures} $\simstruct_i(z)=\simstruct_i^*(z)$ for
all $z\in\defo$.
\end{lem}

\dim{same:structures} Recall first that $\simstruct_i(z)$ is obtained by
glueing together polygons $\tilde\sigma_j(z)$ as in
Lemma~\ref{complex:affine:triangles}.  Moreover, since $\struct(z)$ is
obtained by glueing together the polyhedra $\tilde P_\alpha(z)$ of
Section~\ref{deformation:section}, to get $\simstruct_i^*(z)$ one has to
intersect the $\tilde P_\alpha(z)$ with horospheres centred at ideal
vertices corresponding to the $i$-th cusp, and patch together the
resulting affine polyhedra, which we denote by $Q_l(z)$. 

Both the $\tilde\sigma_j(z)$ and the $Q_l(z)$ are obtained by grouping
together some of the triangles with moduli $z_1,\dots,z_n$, in such a way
that each flat or negative triangle gets grouped with at least one fat
triangle. The grouping rules, however, are different, so indeed we have
something to prove. We first remark that
$\simstruct_i(\zzero)=\simstruct_i^*(\zzero)$, because geometrically (even
if not combinatorially) each $Q_l(\zzero)$ is obtained by glueing together
some $\tilde\sigma_j(\zzero)$.

We will now show that $\simstruct_i(z)$ and $\simstruct_i^*(z)$ have the
same holonomy for $z\in\defo$. Since this holonomy depends analytically on
$\in \defo$, by Lemma~\ref{holonomy}, knowing that
$\simstruct_i(\zzero)=\simstruct_i^*(\zzero)$, it follows from
Theorem~1.7.1 of \cite{CEG} or from \cite{weil} that
$\simstruct_i(z)=\simstruct_i^*(z)$ for $z\in\defo$. Using the recipe
(based on ratios of segments)  mentioned in Lemma~\ref{holonomy}, one gets
combinatorial rules for the holonomies of $\simstruct_i(z)$ and
$\simstruct_i^*(z)$. These rules involve only the moduli $z_1,\dots,z_n$
and apply to loops which are simplicial in the CW-structures on $T_i$
induced respectively by the $\tilde\sigma_j$'s and by the $Q_l$'s. These
CW-structures have, as a common subdivision, the triangulation $\TT_i$
induced by $\TT$ on $T_i$. Using the consistency relations $\cc^*_{\TT_i}$
one easily sees that the two rules extend to one and the same
combinatorial rule which applies to loops which are simplicial in $\TT_i$.
This shows that the holonomies are the same, whence the conclusion.
\finedim{same:structures}

Since $\struct_i(\zzero)$ is a complete cusp, it is isometric to the
quotient of a horoball under the action of $\pi_1(C_i)$ via the holonomy
representation (see Chapter D in \cite{libro} for instance).  Hence, if we
assume that the horoball is centred at
$\infty\in\mc\cup\{\infty\}\cong\partial\mh^3$, the complete cusp has a
developing map of the following form:
$$\begin{array}{rcl} \mr^2\times[0,\infty)&\to& \mh^3\cong\mc\times(0,\infty)\\
(x,y,t)&\mapsto &(a_i(\zzero) x+ b_i(\zzero) y,\exp(t)), \end{array} $$ 
where $a_i(\zzero)$ and $b_i(\zzero)$ are as in Lemma~\ref{holonomy}.

We will apply~\cite{CEG} as in the proof of
Proposition~\ref{develop:torus}. To do this, we shall describe the
holonomy representation of
$\struct_i(z)$ for $z\in\defo$ using the similarity structure on $T_i$
induced by $\struct(z)$. By Lemma~\ref{same:structures}, this structure is
$\simstruct_i(z)$, which is defined as in Proposition~\ref{develop:torus}.
Hence, a holonomy representation for $\struct_i(z)$ can be recovered from
the holonomy representation of $\simstruct_i(z)$ as in
Lemma~\ref{holonomy}, because the hyperbolic holonomy is the conformal
extension of the similarity holonomy.

Then, using~\cite{CEG} as in the proof of Proposition~\ref{develop:torus},
after composing with a hyperbolic isometry we deduce that the following is
a developing map of $\struct_i(z)$ on the cusp $C_i$: 
$$
 \begin{array}{rcl}
 \mr^2\times[0,\infty)&\to&\mh^3\cong\mc\times(0,\infty)\\
 (x,y,t)&\mapsto &
 \left\{\begin{array}{ll}
 ({\exp(u_i(z) x+v_i(z) y)}, \exp(t+\Re ( u_i(z) x+ v_i(z) y)))
  &\hbox{ if } u_i(z)\neq 0\\
 (a_i(z) x+ b_i(z) y,\exp(t) )& \hbox{ otherwise. }
 \end{array}\right.
 \end{array}
$$
Since the argument of~\cite{CEG} applies only to compact manifolds, we
apply it to $\mr^2\times[0,t_n]$ and we consider the limit when
$t_n\to\infty$. This proves the first assertion of the proposition. 

When $u_i(z)=0$, it follows from the expression of this developing map that the
end is complete, as proved in~\cite{libro},~\cite{rat}
or~\cite{thurston:notes}. 

Assume from now to the end of the proof that $u_i(z)\neq 0$. The image of
$\mr^2\times\{t\}$ is precisely the set of points that are at a fixed
distance from the geodesic $\gamma$ with endpoints $0$ and $\infty$, which
is the geodesic fixed by the holonomy representation. Actually, this
distance tends to $0$ as $t$ goes to $\infty$. More precisely, the image
of $\mr^2\times[t,\infty)$ is exactly $U_{r(t)}(\gamma)\setminus\gamma$,
where $U_r$ denotes the tubular $r$-neighbourhood, and $r(t)\to 0$ as
$t\to\infty$.

Let $n_i,m_i\in\mz$ be such that $p_i\cdot n_i-q_i\cdot m_i=1$. The
quadrilateral $Q\subset\mr^2$ with vertices $(0,0)$, $(p_i,q_i)$,
$(p_i+m_i,q_i+n_i)$ and $(m_i,n_i)$ is a fundamental domain for the action
of $\mz^2$ on $\mr^2$. We can also describe $Q$ as: 
$$Q=\{(x,y)\in\mr^2:\ 0\leq n_ix-m_iy\leq 1,\ 0\leq -q_ix+p_iy\leq 1\}.$$
The orbit of $Q$ under the action of the cyclic group 
generated by $(m_i,n_i)$, which corresponds to
$m_i\lambda_i+n_i\mu_i$ in $\pi_1(T_i)$,
is the strip $S=\{(x,y)\in\mr^2:\ 0\leq n_ix-m_iy\leq 1\}$.

First we deal with the case where the relation $p_i\cdot u_i(z)+q_i\cdot
v_i(z)=2\pi\imunit$ is satisfied.  For fixed $t\in[0,\infty)$, the
restriction of the developing map to $S\times\{t\}$ glues one side of $S$
to the other one, and its image is precisely $\partial U_{r(t)}(\gamma)$,
{\it i.e.}\ the set of points at distance $r(t)$ from $\gamma$.  In other
words, the developing map restricted to $\mr^2\times\{t\}$ induces the
universal covering of the cylinder $\partial U_{r(t)}(\gamma)$, and the
deck transformation group is the cyclic group generated by $(p_i,q_i)$,
which corresponds to $p_i \lambda_i + q_i \mu_i$ in $\pi_1(T_i)$.
This description
implies that $C_i$ is isometric to the quotient of
$U_{r(0)}(\gamma)\setminus\gamma$ under the action of the holonomy of
$m_i\lambda_i+n_i\mu_i$. This action extends to a discrete and free action
on the
whole of $U_{r(0)}(\gamma)$, so the completion of $C_i$ is obtained by
adding the quotient of $\gamma$, and the result is a genuine hyperbolic
manifold. Topologically, this manifold is precisely the Dehn filling with
meridian $p_i\lambda_i+q_i\mu_i$. 

In the general case we have $p_i\cdot u_i(z)+q_i\cdot
v_i(z)=2\pi\imunit/r_i$, and we replace $\mh^3$ by a singular space
denoted by $\mh^3_{\alpha_i}$, where $\alpha_i= 2\pi/r_i$.  The space
$\mh^3_{\alpha_i}$ has a singular line $\Sigma\cong\mr$,
$\mh^3_{\alpha_i}\setminus\Sigma$ has a non-complete hyperbolic metric and
the singularity on $\Sigma$ is conical with angle $\alpha_i=2\pi/r_i$.  In
cylindrical coordinates the metric on $\mh^3_{\alpha_i}\setminus\Sigma$
has the form:
$$ {\rm d}s^2={\rm d}r^2+\left({\alpha_i\over 2\pi}\right)^2\sinh^2 (r) 
{\rm d}\vartheta^2+ \cosh^2(r) {\rm d}h^2 $$
where $r\in (0,+\infty)$ is the distance to $\Sigma$, $\vartheta\in
[0,2\pi)$ is the angular parameter and $h\in\mr$ if the height.

The developing
map $\widetilde C_i\to \mh^3\setminus\gamma$ induces a developing map
$\widetilde C_i\to \mh^3_{\alpha_i}\setminus\Sigma$, because the universal
coverings of ${\mh^3\setminus\gamma}$ and of
${\mh^3_{\alpha_i}\setminus\Sigma}$ are isometric. Then the argument in
the non-singular case above (where $r_i=1$) applies to the singular case
after replacing the pair $(\mh^3,\gamma)$ by $(\mh^3_{\alpha_i},\Sigma)$.
The completion is of course in this case a cone manifold with cone angle
$\alpha_i$ along the loop added.\finedim{develop:cusp}

Proposition~\ref{develop:cusp} and the discussion preceding it imply
point~\ref{completion:description:point} in Theorem~\ref{completion:teo}.
We are only left to establish point~\ref{rigidity:u:point}, which we do now.

\begin{prop}\label{rigidity:prop}
If $z\in\defo$ and $u_1(z)=\cdots=u_k(z)=0$ then $z=\zzero$.
\end{prop}

\dim{rigidity:prop} Having already established point~\ref{symmetric:u:v:point}
in Theorem~\ref{completion:teo}, we can rephrase the statement as follows: {\em
if $\defo_0$ is the set of solutions $z$ in $\zznbd$ of both $\cc_\TT$ and
$\mm_\TT$, then $\zzero$ is an isolated point of $\defo_0$.} Assume this is not
the case. Since $\defo_0$ is an analytic space, we can find a non-constant
curve in $\defo_0$ starting at $\zzero$. Therefore, at least one of the
coordinates $z_j$ assumes uncountably many different values on $\defo_0$.

Now, by Proposition~\ref{develop:cusp}, every $z\in\defo_0$ defines on $M$ a
complete finite-volume hyperbolic structure, which must be isometric to the
original structure by Mostow rigidity. It follows that for all $z\in\defo_0$
the original manifold $M$ contains a geodesic ideal tetrahedron, possibly flat
and with some paired faces, of modulus $z_j$ or $\overline{z_j}$, depending on
whether $\Im(z_j)$ is non-negative or non-positive. In particular, under
the
assumption that $\zzero$ is not isolated, $M$ contains uncountably many
pairwise non-isometric (possibly flat) geodesic ideal tetrahedra.

Let us consider now the universal covering $\mh^3\to M$, on which the group of
deck transformations acts as a subgroup of $\isompiutre$ identified to
$\pi_1(M)$. It is very easy to see that each geodesic ideal tetrahedron
contained in $M$ is actually the projection of the convex hull of 4 points of
$\partial\mh^3$ which are fixed points of parabolic elements of $\pi_1(M)$.
Since $\pi_1(M)$ is countable and each parabolic element has one fixed point,
we see that in $M$ there are at most countably many pairwise non-isometric
(possibly flat) geodesic ideal tetrahedra. This gives a contradiction and
concludes the proof.\finedim{rigidity:prop}

\section{Hyperbolic filling parameters}\label{parameters:section}

The aim of this section is to show that the set of parameters $(c_1,\dots,c_k)$
arising as in
Theorem~\ref{completion:teo}(\ref{completion:description:point}) covers a
neighbourhood of $(\infty,\dots,\infty)$ in $(\mz^2\sqcup\{\infty\})^k$. This
will imply the conclusion of the proof. We will start with a combinatorial
argument due to  Neumann and Zagier~\cite{ne:za}, which shows that the space
$\defo$ of deformed structures is sufficiently big ({\it i.e.}~it has
(complex) dimension exactly $k$). Later we will modify the approach
of~\cite{ne:za} to
avoid the assumption that $\zzero$ is a smooth point of $\defo$.

Note first that the expressions $z$, $1/(1-z)$ and $1-1/z$ can all be rewritten
as $\delta_0\cdot z^{\delta_1}\cdot (1-z)^{\delta_2}$ for suitable $\delta_0,
\delta_1,\delta_2\in\{\pm 1\}$. Recall that our ideal triangulation $\TT$ of
$M$ consists of tetrahedra $\Delta_j$, $j=1,\dots,n$, and $\partial M$ consists
of tori $T_i$, $i=1,\dots,k$.

\begin{lem}\label{chi:counts:edges}
$\TT$ contains $n$ edges.
\end{lem}

\dim{chi:counts:edges} Since $\partial\Mbar$ is made of tori, $\chi(\Mbar)=0$.
Hence $\chi(\hatM)=k$, because each torus is collapsed to a point. In $\TT$
there are twice as many faces as tetrahedra, so $k=k-(\#{\rm edges})+2n-n$,
whence the conclusion.\finedim{chi:counts:edges}

Let us list the edges in $\TT$ as $e_m$, $m=1,\dots,n$. For
$m,j\in\{1,\dots,n\}$ let us define $(\theta_1(m,j),\theta_2(m,j))$ as the sum
of all pairs $(\delta_1,\delta_2)$ over the edges $e$ of $\Delta_j$ which get
identified to $e_m$, where the modulus of $\Delta_j$ along $e$ is  $\pm
z_j^{\delta_1}(1-z_j)^{\delta_2}$. For suitable $\varepsilon_m\in\{\pm 1\}$,
$m=1,\dots,n$, we can therefore write $\cc_\TT(z)$ as 
\begin{equation}
\prod_{j=1}^n z_j^{\theta_1(m,j)}\cdot(1-z_j)^{\theta_2(m,j)}=\varepsilon_m,
\qquad m=1,\dots,n.\label{theta:system}
\end{equation}
Let us denote now by $v_i$ the vertex of $\hatM$ obtained by collapsing
$T_i\subset\partial\Mbar$. For $i\in\{1,\dots,k\}$ and $m\in\{1,\dots,n\}$ we
define $x(i,m)\in\{0,1,2\}$ as the number of ends of $e_m$ which get identified
to $v_i$ in $\Mbar$. We have now two matrices $X\in\mm(k\times n,\mc)$ and 
$\Theta=(\Theta_1,\Theta_2)\in\mm(n\times 2n,\mc)$. The entries are actually
integers, but it will be convenient to view $X$ and $\Theta$ as complex
matrices. The next two combinatorial results are due to Neumann and
Zagier~\cite{ne:za} and show that $\defo$ is an open portion of a complex 
algebraic
variety of dimension at least $k$. We note that this result in~\cite{libro} was
deduced from a much harder combinatorial lemma from~\cite{ne:za}.

\begin{lem}\label{XTheta=0} $X\cdot\Theta=0$.
\end{lem}

\dim{XTheta=0} We must check that for all $i$ and $j$
\begin{eqnarray*}
& & \sum_{m=1}^n x(i,m)\cdot\theta_1(m,j)=
\sum_{m=1}^n x(i,m)\cdot\theta_2(m,j)=0\\
&{\it i.e.}\quad &\sum_{m=1}^n x(i,m)\cdot(\theta_1(m,j),\theta_2(m,j))=0.
\end{eqnarray*}
We can rewrite the last sum as
\begin{eqnarray*}
& & 	\sum_{m=1}^n \quad
	\sum_{
		{\tiny\begin{array}{c} 
		v{\rm\ endpoint\ of\ }e_m\\ 
		v{\rm\ identified\ to\ }v_i
		\end{array}}}\quad
	\sum_{
		{\tiny\begin{array}{c} 
		e{\rm\ edge\ of\ }\Delta_j\\
		e{\rm\ identified\ to\ }e_m\\
		{\rm mod}(\Delta_j|e)=\pm z_j^{\delta_1}(1-z_j)^{\delta_2}
		\end{array}}}\quad
	(\delta_1,\delta_2)\\
& = & 	\sum_{
		{\tiny\begin{array}{c}
		v{\rm\ vertex\ of\ }\Delta_j\\
		v{\rm\ identified\ to\ }v_i
		\end{array}}}\quad
	\sum_{
		{\tiny\begin{array}{c} 
		e{\rm\ edge\ of\ }\Delta_j\\
		e{\rm\ contains\ }v{\rm\ as\ endpoint}\\
		{\rm mod}(\Delta_j|e)=\pm z_j^{\delta_1}(1-z_j)^{\delta_2}
		\end{array}}}\quad
	(\delta_1,\delta_2)\\
& = & 	\sum_{
		{\tiny\begin{array}{c}
		v{\rm\ vertex\ of\ }\Delta_j\\
		v{\rm\ identified\ to\ }v_i
		\end{array}}}\quad
	\Big((1,0)+(0,-1)+(-1,1)\Big)=0.
\end{eqnarray*}
This concludes the proof.\finedim{XTheta=0}

\begin{lem}\label{rankX=k}{\rm rank}$_{\,\mc}(X)=k$.
\end{lem}

\dim{rankX=k}
Let $a_1,\dots,a_k\in\mc$ be such that $(a_1,\dots,a_k)\cdot X=0$, {\em i.e.}
$$\sum_{i=1}^ka_i\cdot x(i,m) =0,\qquad m=1,\dots,n.$$
Using the definition of $X$, this means that $a_{i_0}+a_{i_1}=0$ whenever 
$v_{i_0}$ and $v_{i_1}$ are the ends of some edge in $\hatM$.
If we examine a face of some $\Delta_j$ having vertices $v_{i_0}$, $v_{i_1}$ and
$v_{i_2}$, the three edges of the face yield respectively the
relations 
$$a_{i_0}+a_{i_1}=0,\qquad a_{i_0}+a_{i_2}=0,
\qquad a_{i_1}+a_{i_2}=0.$$
Therefore $a_i=0$ for $i=1,\dots,k$, and the conclusion follows.
\finedim{rankX=k}

\begin{cor}{\rm rank$_{\,\mc}(\Theta)\leq n-k$, in particular $k\leq n$.}
\end{cor}

Going back to the system $\cc_\TT$ written as in formula~(\ref{theta:system}),
we can now show that it can be replaced by a system of $n-k$ equations only.
This fact, even if not explicitly stated in~\cite{ne:za}, was certainly known
to the authors. We reproduce here with minor improvements the proof given
in~\cite{libro}. For the sake of simplicity we rearrange the edges
$e_1,\dots,e_n$ in such a way that the last $k$ rows of $\Theta$ are linearly
dependent on the first $n-k$.

\begin{prop}\label{defo:n-k:equations}
$$\defo=\left\{z\in\uu: 
\prod_{j=1}^n z_j^{\theta_1(m,j)}\cdot(1-z_j)^{\theta_2(m,j)}=\varepsilon_m,
\ m=1,\dots,n-k\right\}.$$
\end{prop}

\dim{defo:n-k:equations}
We can choose continuous branches of the logarithm function near 
$\zzero_j$ and $(1-\zzero_j)$, $j=1,\dots,n$, and assume that the neighbourhood
$\zznbd$ of $\zzero$ used to define $\defo$ is small
enough that $\log(z_j)$ and $\log(1-z_j)$ are defined for $z\in\zznbd$.
By the properties of the exponential map there exist constants
$r_m\in\mz$, $m=1,\dots,n$, such that
$$\sum_{j=1}^n\Big(\theta_1(m,j)\log(\zzero_j)+
\theta_2(m,j)\log(1-\zzero_j)\Big)=\imunit\pi(2r_m+(1-\varepsilon_m)/2).$$
By continuity, if $\zznbd$ is small enough, for $z\in\zznbd$ and 
$m\in\{1,\dots,n\}$ the next two equations
are equivalent:
\begin{eqnarray}
\prod_{j=1}^n z_j^{\theta_1(m,j)}\cdot(1-z_j)^{\theta_2(m,j)}&=&\varepsilon_m,\\
\sum_{j=1}^n\Big(\theta_1(m,j)\log(z_j)+
\theta_2(m,j)\log(1-z_j)\Big)&=&
\imunit\pi(2r_m+(1-\varepsilon_m)/2).\label{log:form:eqns}
\end{eqnarray}
We have to show that the first $n-k$ of these equations imply the last $k$ of
them. We will use the logarithm form~(\ref{log:form:eqns})
 of the equations. By assumption,
for $m>n-k$ there exist $a^1_m,\dots,a^{n-k}_m\in\mc$ such that
$$\theta_t(m,j)=\sum_{l=1}^{n-k} a^l_m\cdot\theta_t(l,j),\qquad t=1,2,
\quad j=1,\dots,n.$$
Therefore if $z\in\zznbd$ solves the first $n-k$ equations we have for $m>n-k$
\begin{eqnarray*}
& & \sum_{j=1}^n\Big(\theta_1(m,j)\log(z_j)+
\theta_2(m,j)\log(1-z_j)\Big)\\
& = & \sum_{j=1}^n\sum_{l=1}^{n-k} a^l_m
\Big(\theta_1(l,j)\log(z_j)+
\theta_2(l,j)\log(1-z_j)\Big)\\
& = & \sum_{l=1}^{n-k} a^l_m
\imunit\pi(2r_l+(1-\varepsilon_l)/2).
\end{eqnarray*}
For $z=\zzero$ the first line equals $\imunit\pi(2r_m+(1-\varepsilon_m)/2)$
so the last line has the same (constant) value, and the conclusion follows.
\finedim{defo:n-k:equations}

We note now that by
Theorem~\ref{completion:teo}(\ref{symmetric:u:v:point},\ref{non:real:limit:point})
for $i=1,\dots,k$ we can define a function
$g_i:\defo\to S^2=\mr^2\sqcup\infty$
as $g_i(z)=\infty$ if $u_i(z)=0$, and $g_i(z)$ as the only pair $(p,q)$ of real
numbers such that $p\cdot u_i(z)+q\cdot v_i(z)=2\pi\imunit$ otherwise. The rest
of this section is devoted to establishing the following result, which,
together with
Theorem~\ref{completion:teo}(\ref{completion:description:point}), implies
Theorem~\ref{new:main:teo} and hence Theorem~\ref{main:teo}.

\begin{prop}\label{infinity:covered:prop}
The image of $g=(g_1,\dots,g_k):\defo\to(S^2)^k$ covers a neighbourhood of 
$(\infty,\dots,\infty)$.
\end{prop}

\dim{infinity:covered:prop} Let us consider the homeomorphism $\varphi_i:S^2\to
S^2$ defined by  $$\varphi_i(p,q)={2\pi\imunit\over p+\tau_i q}$$ (we are
viewing the first $S^2$ as $\mr^2\sqcup\{\infty\}$ and the second one as
$\mc\sqcup\{\infty\}$, and as usual $1/0=\infty$, $1/\infty=0$). We define now
$\tilde u_i:\defo\to\mc$ as $\varphi_i\compo g_i$. To conclude it is sufficient
to show that the image of $\tilde u=(\tilde u_1,\dots,\tilde
u_k):\defo\to\mc^k$ covers a neighbourhood of $0$.
Recall first the following two essential properties of $\defo$ already
established:
\begin{enumerate}
\item $\defo$ is a (germ of) analytic variety, defined in $\mc^n$ as the zero
set of $n-k$ holomorphic functions.
\item There is a map $u:\defo\to\mc^k$ which is the restriction of a holomorphic
function on an open subset of $\mc^n$, such that $u^{-1}(\{0\})=\{\zzero\}$.
\end{enumerate}
Under these assumptions, the preparation theorem of
Weierstrass~\cite{milnor:sing} implies that $u:\defo\to\mc^k$ is an open map
(more precisely, it is a covering branched over a {\em real} codimension-2
set).
We denote now by $\|\,.\,\|$ the usual Euclidean norm on $\mc^k$, and claim
that 
\begin{equation}
\lim_{z\in\defo,z\to\zzero}{\|\tilde u-u\| \over \|u\|}=0.\label{zero:limit}
\end{equation}
Of course it is sufficient to show that for all $i$
$$\lim_{z\in\defo,u_i(z)\neq 0,z\to\zzero}{\tilde u_i-u_i \over u_i}=0.$$
Using the relations
$$p_i\cdot u_i+q_i\cdot v_i=2\pi\imunit,
\qquad\tilde u_i={2\pi\imunit\over p_i+\tau_iq_i},
\qquad{v_i\over u_i}\longrightarrow\tau_i,
\qquad \Im(\tau_i)\neq0$$
we see that 
$${\tilde u_i-u_i \over u_i}=\tau_i\cdot {q_i\over p_i+\tau_iq_i}\cdot
\left({v_i\over\tau_iu_i}-1\right)\longrightarrow 0$$
because $|q_i/(p_i+\tau_iq_i)|$ is bounded from above by $1/|\Im(\tau_i)|$.
Formula~(\ref{zero:limit}) is proved.

Let us consider now the function $\|u\|:\defo\to\mr_+$, denoted by $f$. Note
that $f^{-1}(\{0\})=\{\zzero\}$ and that $f$ is the restriction to $\defo$ of
an ambient map whose square is real-analytic. Since $u:\defo\to\mc^k$ is open
and $u^{-1}(\{0\})=\{\zzero\}$, we can choose a small $R>0$ and restrict
$\defo$ so that $u:\defo\to B_R(0)$ is proper and surjective. Here $B_R(0)$ is
the open $R$-ball centred at $0$ in $\mc^k$. In the sequel $S_r(0)$ will denote
the $R$-sphere. For $0<r<R$ we also set $\defo_{\leq r}=f^{-1}([0,r])$, and
$\defo_{=r}=f^{-1}(\{r\})$.

Using the general theory of analytic spaces~\cite{whitney} we can now choose a
good stratification of $\defo$ (with respect to singularity), and assume that,
away from $\zzero$, $f$ is transversal to all strata. Since $\defo$ is defined
by complex-analytic functions, its top real $2k$-dimensional strata are
naturally
oriented, and there are no strata of real dimension $2k-1$. Now, for
$0<r<R$ we
consider the induced stratification of $\defo_{=r}$, and orient the top
real $(2k-1)$-dimensional strata using $f$ and the previous orientation.
($\defo_{\leq r}$ actually has a (stratified) conic structure with basis
$\defo_{=r}$, vertex $\zzero$ and height function $f$, but we will not need all
this information.) Since in $\defo_{=r}$ there are no strata of
real dimension
$2k-2$, we can view it as a geometric $(2k-1)$-cycle. Similarly, $\defo_{\leq
r}$ can be regarded as a $2k$-dimensional geometric $\mz$-chain with boundary
$\defo_{=r}$. 

Using Sard's lemma we see that there exist arbitrarily small $r>0$ and regular
values $w\in S_r(0)$ for the restriction of $u$ to all strata of both $\defo$
and $\defo_{=r}$. Since $u$ is complex-analytic, each preimage of $w$ has index
$+1$ with respect to $u$. Orientation conventions imply that the same is true
with respect to $u|_{\defo_{=r}}$. Moreover $u|_{\defo_{=r}}$ is surjective
onto $S_r(0)$. We deduce that $u|_{\defo_{=r}}:\defo_{=r}\to S_r(0)$, as a
geometric cycle in $S_r(0)$, represents a {\em strictly positive} (in
particular, non-zero) multiple of the canonical generator of
$H_{2k-1}(S_r(0))\cong H_{2k-1}(\mc^k\setminus\{0\})$ (we will take all
homology groups with integer coefficients).

Using formula~(\ref{zero:limit}) we can now assume that $r$ is small enough
that $\|\tilde u-u\|<\|u\|/2$ on $\defo_{\leq r}$, in particular  $\|\tilde
u-u\|< r/2$.  This implies that $\tilde
u(\defo_{=r})\subset\mc^k\setminus\{0\}$, moreover $\tilde
u:\defo_{=r}\to\mc^k\setminus\{0\}$ is homotopic to
$u:\defo_{=r}\to\mc^k\setminus\{0\}$, whence it represents the same non-zero
element of $H_{2k-1}(\mc^k\setminus\{0\})$.

We claim now that $\tilde u(\defo_{\leq r})$ contains $D_{r/2}(0)$.  Assume by
contradiction that there exists $w^{(0)}\in D_{r/2}(0)\setminus\tilde
u(\defo_{\leq r})$.  Note that each half-line in $\mc^k$ with origin in
$w^{(0)}$ meets $S_r(0)$ exactly once, so we have a natural ``radial''
projection $p:\mc^k\setminus\{w^{(0)}\}\to S_r(0)$. We can now consider the
$2k$-dimensional geometric chain $p\compo\tilde u:\defo_{\leq r}\to S_r(0)$,
whose boundary $p\compo\tilde u:\defo_{=r}\to S_r(0)$ is therefore zero in
$H_{2k-1}(S_r(0))\cong H_{2k-1}(\mc^k\setminus\{0\})$. Now, for
$z\in\defo_{=r}$ we have $\|\tilde u(z)\|>r/2$. Since $\|w^{(0)}\|\leq r/2$, by
the definition of $p$, the segment joining $\tilde u(z)$ and $p(\tilde u(z))$
does not contain $0$. In particular, the geometric $(2k-1)$-chains
$p\compo\tilde u:\defo_{=r}\to\mc^k\setminus\{0\}$ and $\tilde
u:\defo_{=r}\to\mc^k\setminus\{0\}$ are homotopic to each other in
$\mc^k\setminus\{0\}$. This is a contradiction, because the former is zero in
$H_{2k-1}(\mc^k\setminus\{0\})$ and the latter is not. Our claim is
established and the proof is complete. \finedim{infinity:covered:prop}

\vspace{1.5cm}

\hspace{8cm} Dipartimento di Matematica
\vspace{-1pt}

\hspace{8cm} Universit\`a di Pisa
\vspace{-1pt}

\hspace{8cm} Via F. Buonarroti, 2
\vspace{-1pt}

\hspace{8cm} I-56127 Pisa, Italy
\vspace{-1pt}

\hspace{8cm} petronio@dm.unipi.it

\vspace{.5cm}

\hspace{8cm} 	Universitat Aut\`onoma de Barcelona
\vspace{-1pt}

\hspace{8cm} 	Departament de Matem\`atiques
\vspace{-1pt}

\hspace{8cm} 	E-08193 Bellaterra, Spain
\vspace{-1pt}

\hspace{8cm} porti@manwe.mat.uab.es

\end{document}